\documentclass[envcountsame]{llncs}

\usepackage{amsmath}
\usepackage{amssymb}

\newcommand{\psod}[1]{~({#1})}

\newcommand{\C}{\mathbb C}

\newcommand{\Q}{\mathbb Q}
\newcommand{\R}{\mathbb R}
\newcommand{\Z}{\mathbb Z}

\newcommand{\eps}{\epsilon}

\DeclareMathOperator{\disc}{disc}

\DeclareMathOperator{\Tr}{Tr}

\begin{document}

\title{Enumeration of totally real number fields \\ of bounded root discriminant}
\author{John Voight}
\institute{Department of Mathematics and Statistics \\ University of Vermont \\ Burlington, VT\ 05401 \\ \email{jvoight@gmail.com}}

\maketitle

\begin{abstract}
We enumerate all totally real number fields $F$ with root discriminant $\delta_F \leq 14$.  There are $1229$ such fields, each with degree $[F:\Q] \leq 9$.
\end{abstract}

In this article, we consider the following problem.

\begin{problem}
Given $B \in \R_{>0}$, enumerate the set $NF(B)$ of totally real number fields $F$ with root discriminant $\delta_F\leq B$, up to isomorphism.
\end{problem}

To solve Problem 1, for each $n \in \Z_{>0}$ we enumerate the set
\[ NF(n,B)=\{F \in NF(B):[F:\Q]=n\} \]
which is finite (a result originally due to Minkowski).  If $F$ is a totally real field of degree $n=[F:\Q]$, then by the Odlyzko bounds \cite{Odlyzko}, we have $\delta_F \geq 4\pi e^{1+\gamma} - O(n^{-2/3})$ where $\gamma$ is Euler's constant; thus for $B < 4\pi e^{1+\gamma} < 60.840$, we have $NF(n,B)=\emptyset$ for $n$ sufficiently large and so the set $NF(B)$ is finite.  Assuming the generalized Riemann hypothesis (GRH), we have the improvement $\delta_F \geq 8 \pi e^{\gamma+\pi/2} - O(\log^{-2} n)$ and hence $NF(B)$ is conjecturally finite for all $B < 8\pi e^{\gamma+\pi/2} < 215.333$.  On the other hand, for $B$ sufficiently large, the set $NF(B)$ is infinite: Martin \cite{Martin} has constructed an infinite tower of totally real fields with root discriminant $\delta \approx 913.493$ (a long-standing previous record was held by Martinet \cite{Martinet1} with $\delta \approx 1058.56$).  The value
\[ \liminf_{n \to \infty}\, \min \{\delta_F:F\in NF(n,B)\} \]
is presently unknown.  If $B$ is such that $\#NF(B)=\infty$, then to solve Problem 1 we enumerate the set $NF(B)=\bigcup_n NF(n,B)$ by increasing degree. 

Our restriction to the case of totally real fields is not necessary: one may place alternative constraints on the signature of the fields $F$ under consideration (or even analogous $p$-adic conditions).  However, we believe that Problem 1 remains one of particular interest.  First of all, it is a natural boundary case: by comparison, Hajir-Maire \cite{HajirMaire,HajirMaire2} have constructed an unramified tower of totally complex number fields with root discriminant $\approx 82.100$, which comes within a factor $2$ of the GRH-conditional Odlyzko bound of $8\pi e^{\gamma} \approx 44.763$.  Secondly, in studying certain problems in arithmetic geometry and number theory---for example, in the enumeration of arithmetic Fuchsian groups \cite{LMR} and the computational investigation of the Stark conjecture and its generalizations---provably complete and extensive tables of totally real fields are useful, if not outright essential.  Indeed, existing strategies for finding towers with small root discriminant as above often start by finding a good candidate base field selected from existing tables.

The main result of this note is the following theorem, which solves Problem 1 for $\delta=14$.

\begin{theorem}
We have $\#NF(14)=1229$.
\end{theorem}

The complete list of these fields is available online \cite{VoightNumberfields}; the octic and nonic fields ($n=8,9$) are recorded in Tables 3--4 in \S 4, and there are no dectic fields ($NF(14,10)=\emptyset$).  For a comparison of this theorem with existing results, see \S 1.2.

The note is organized as follows.  In \S 1, we set up the notation and background.  In \S 2, we describe the computation of primitive fields $F \in NF(14)$; we compare well-known methods and provide some improvements.  In \S 3, we discuss the extension of these ideas to imprimitive fields, and we report timing details on the computation.  Finally, in \S 4 we tabulate the fields $F$.

The author wishes to thank: J\"urgen Kl\"uners, Noam Elkies, Claus Fieker, Kiran Kedlaya, Gunter Malle, and David Dummit for useful discussions; William Stein, Robert Bradshaw, Craig Citro, Yi Qiang, and the rest of the Sage development team for computational support (NSF Grant No.~0555776); and Larry Page and Helen Read for their technical assistance.  

\section{Background}

\subsection{Initial bounds}

Let $F$ denote a totally real field of degree $n=[F:\Q]$ with discriminant $d_F$ and root discriminant $\delta_F=d_F^{1/n}$.  By the unconditional Odlyzko bounds \cite{Odlyzko} (see also Martinet \cite{Martinet}), if $n \geq 11$ then $\delta_F > 14.083$, thus if $F \in NF(14)$ then $n \leq 10$.  

\begin{table}[h] 
\caption{Degree and Root Discriminant Bounds}
\begin{center}
\begin{tabular}{c|ccccccccc}
$n$ & 2 & 3 & 4 & 5 & 6 & 7 & 8 & 9 & 10 \\
\hline
$B_O$ & $>$ 2.223 & 3.610 & 5.067 & 6.523 & 7.941 & 9.301 & 10.596 & 11.823 & 12.985 \\
$B_O$ (GRH) & $>$ 2.227 & 3.633 & 5.127 & 6.644 & 8.148 & 9.617 & 11.042 & 12.418 & 13.736 \\
\hline
$\Delta$ & 30 & 25 & 20 & 17 & 16 & 15.5 & 15 & 14.5 & 14
\end{tabular}
\end{center}
\end{table}

The lower bounds for $\delta_F$ in the remaining degrees are summarized in Table 1: for each degree $2 \leq n \leq 10$, we list the unconditional Odlyzko bound $B_O=B_O(n)$, the GRH-conditional Odlyzko bound (for comparison only, as computed by Cohen-Diaz y Diaz-Olivier \cite{CDO}), and the bound $\delta_F \leq \Delta$ that we employ.

\subsection{Previous work}

There has been an extensive amount of work done on the problem of enumerating number fields---we refer to \cite{KMart} for a discussion and bibliography.

\begin{enumerate}
\item The KASH and PARI groups \cite{KASH} have computed tables of number fields of all signatures with degrees $\leq 7$: in degrees $6,7$, they enumerate totally real fields up to discriminants $10^7,15\cdot 10^7$, respectively (corresponding to root discriminants $14.67,14.71$, respectively).  
\item Malle \cite{Malle109} has computed all totally real primitive number fields of discriminant $d_F \leq 10^9$ (giving root discriminants $31.6, 19.3, 13.3, 10$ for degrees $6,7,8,9$).  This was reported to take several years of CPU-time on a SUN workstation.
\item The database by Kl\"uners-Malle \cite{KM} contains polynomials for all transitive groups up to degree $15$ (including possible combinations of signature and Galois group); up to degree 7, the fields with minimal (absolute) discriminant with given Galois group and signature have been included.
\item Roblot \cite{Roblot} constructs abelian extensions of totally real fields of degrees $4$ to $48$ (following Cohen-Diaz y Diaz-Olivier \cite{CDyD}) with small root discriminant.
\end{enumerate}

The first two of these allow us only to determine $NF(10)$ (if we also separately compute the imprimitive fields); the latter two, though very valuable for certain applications, are in a different spirit than our approach.  Therefore our theorem substantially extends the complete list of fields in degrees $7$--$9$.

\section{Enumeration of totally real fields}

\subsection{General methods}

The general method for enumerating number fields is well-known (see Cohen \cite[\S 9.3]{Cohen2}).  We define the Minkowski norm on a number field $F$ by $T_2(\alpha)=\sum_{i=1}^n |\alpha_i|^2$ for $\alpha \in F$, where $\alpha_1,\alpha_2,\dots,\alpha_n$ are the conjugates of $\alpha$ in $\C$.  The norm $T_2$ gives $\Z_F$ the structure of a lattice of rank $n$.  In this lattice, the element $1$ is a shortest vector, and an application of the geometry of numbers to the quotient lattice $\Z_F/\Z$ yields the following result.

\begin{lemma}[Hunter] \label{Hunter}
There exists $\alpha \in \Z_F \setminus \Z$ such that $0 \leq \Tr(\alpha) \leq n/2$ and
\[ T_2(\alpha) \leq \frac{\Tr(\alpha)^2}{n} + \gamma_{n-1}\left(\frac{|d_F|}{n}\right)^{1/(n-1)} \]
where $\gamma_{n-1}$ is the $(n-1)$th Hermite constant.
\end{lemma}

\begin{remark}
The values of the Hermite constant are known for $n \leq 8$ (given by the lattices $A_1,A_2,A_3,D_4,D_5,E_6,E_7,E_8$): we have $\gamma_n^n=1,4/3,2,4,8,64/3,64,256$ (see Conway and Sloane \cite{CS}) for $n=1,\dots,8$; the best known upper bounds for $n=9,10$ are given by Cohn and Elkies \cite{CE}.
\end{remark}

Therefore, if we want to enumerate all number fields $F$ of degree $n$ with $|d_F| \leq B$, an application of Lemma \ref{Hunter} yields $\alpha \in \Z_F \setminus \Z$ such that $T_2(\alpha) \leq C$ for some $C \in \R_{>0}$ depending only on $n,B$.  We thus obtain bounds on the power sums
\[ |S_k(\alpha)|=\biggl|\sum_{i=1}^n \alpha_i^k\biggr| \leq T_k(\alpha)=\sum_{i=1}^n |\alpha_i|^k \leq nC^{k/2}, \]
and hence bounds on the coefficients $a_i \in \Z$ of the characteristic polynomial 
\[ f(x)=\prod_{i=1}^n (x-\alpha_i) = x^n + a_{n-1} x^{n-1} + \dots + a_0 \] 
of $\alpha$ by Newton's relations:
\begin{equation}
S_k+\sum_{i=1}^{k-1}a_{n-1}S_{k-i} + ka_{n-k}=0.
\label{Newton}
\end{equation}
This then yields a finite set $NS(n,B)$ of polynomials $f(x) \in \Z[x]$ such that every $F$ is represented as $\Q[x]/(f(x))$ for some $f(x) \in NS(n,B)$, and in principle each $f(x)$ can then be checked individually.  We note that it is possible that $\alpha$ as given by Hunter's theorem may only generate a subfield $\Q\subset \Q(\alpha) \subsetneq F$ if $F$ is imprimitive: for a treatment of this case, see \S 3.

The size of the set $NS(n,B)$ is $O(B^{n(n+2)/4})$ (see Cohen \cite[\S 9.4]{Cohen2}), and the exponential factor in $n$ makes this direct method impractical for large $n$ or $B$.  Note, however, that it is sharp for $n=2$: we have $NF(2,B) \sim (6/\pi^2)B^2$ (as $B \to \infty$), and indeed, in this case one can reduce to simply listing squarefree integers.  For other small values of $n$, better algorithms are known: following Davenport-Heilbronn, Belabas \cite{Belabas} has given an algorithm for cubic fields; Cohen-Diaz y Diaz-Olivier \cite{CDO} use Kummer theory for quartic fields; and by work of Bhargava \cite{Bhargava}, in principle one should similarly be able to treat the case of quintic fields.  No known method improves on this asymptotic complexity for general $n$, though some possible progress has been made by Ellenberg-Venkatesh \cite{EV}.

\subsection{Improved methods for totally real fields}

We now restrict to the case that $F$ is totally real.  Several methods can then be employed to improve the bounds given above---although we only improve on the implied constant in the size of the set $NS(n,B)$ of examined polynomials, these improvements are essential for practical computations.  

\subsubsection*{Basic bounds.}

From Lemma \ref{Hunter}, we have $0 \leq a_{n-1} = -\Tr(\alpha) \leq \left\lfloor n/2 \right\rfloor$ and 
\[ a_{n-2}=\frac{1}{2}a_{n-1}^2-\frac{1}{2}T_2(\alpha) \geq \frac{1}{2}\left(1-\frac{1}{n}\right) a_{n-1}^2 - \frac{\gamma_{n-1}}{2}\left(\frac{B}{n}\right)^{1/(n-1)}. \]

For an upper bound on $a_{n-2}$, we apply the following result.

\begin{lemma}[{Smyth \cite{Smyth}}] \label{Smyth}
If $\gamma$ is a totally positive algebraic integer, then
\[ \Tr(\gamma)>1.7719[\Q(\gamma):\Q] \]
unless $\gamma$ is a root of one of the following polynomials:
\[ x-1, x^2-3x+1, x^3-5x^2+6x-1, x^4-7x^3+13x^2-7x+1, x^4-7x^3+14x^2-8x+1. \]
\end{lemma}

\begin{remark}
The best known bound of the above sort is due to Aguirre-Bilbao-Peral \cite{ABP}, who give $\Tr(\gamma) > 1.780022[\Q(\gamma):\Q]$ with $14$ possible explicit exceptions.  For our purposes (and for simplicity), the result of Smyth will suffice.
\end{remark}

Excluding these finitely many cases, we apply Lemma \ref{Smyth} to the totally positive algebraic integer $\alpha^2$, using the fact that $T_2(\alpha)=\Tr(\alpha^2)$, to obtain the upper bound $a_{n-2} < a_{n-1}^2/2 - 0.88595n$.

\subsubsection*{Rolle's theorem.}

Now, given values $a_{n-1},a_{n-2},\dots,a_{n-k}$ for the coefficients of $f(x)$ for some $k \geq 2$, we deduce bounds for $a_{n-k-1}$ using Rolle's theorem---this elementary idea can already be found in Takeuchi \cite{Takeuchi} and Kl\"uners-Malle \cite[\S 3.1]{KMart}.  Let
\[ f_i(x)=\frac{f^{(n-i)}(x)}{(n-i)!}=g_i(x) + a_{n-i} \]
for $i=0,\dots,n$.  Consider first the case $k=2$.  Then
\[ g_3(x)=\frac{n(n-1)(n-2)}{6}x^3 + \frac{(n-1)(n-2)}{2} a_{n-1} x^2 + (n-2)a_{n-2}x. \]
Let $\beta_1<\beta_2$ denote the roots of $f_2(x)$.  Then by Rolle's theorem, 
\[ f_3(\beta_1)=g_3(\beta_1)+a_{n-3}>0 \quad \text{and} \quad f_3(\beta_2)=g_3(\beta_2)+a_{n-3}<0 \]
hence $-g_3(\beta_1)<a_{n-3}<-g_3(\beta_2)$.  In a similar way, if $\beta_1^{(k)}<\dots<\beta_k^{(k)}$ denote the roots of $f_k(x)$, then we find that
\[ -\min_{\substack{1 \leq i \leq k \\ i \not\equiv k \psod{2}}} g_{k+1}(\beta_i^{(k)})<a_{n-k-1}< - \max_{\substack{1 \leq i \leq k \\ i \equiv k \psod{2}}} g_{k+1}(\beta_i^{(k)}). \]

\subsubsection*{Lagrange multipliers.}

We can obtain further bounds as follows.  We note that if the roots of $f$ are bounded below by $\beta_0^{(k)}$  (resp.\ bounded above by $\beta_{k+1}^{(k)}$), then 
\[ f_k(\beta_0^{(k)})=g_k(\beta_0^{(k)})+a_{n-k}> 0 \]
(with a similar inequality for $\beta_{k+1}^{(k)}$), and these combine with the above to yield
\begin{equation} \label{Rolle}
-\min_{\substack{0 \leq i \leq k+1 \\ i \not\equiv k \psod{2}}} g_{k+1}(\beta_i^{(k)})<a_{n-k-1}< - \max_{\substack{0 \leq i \leq k+1 \\ i \equiv k \psod{2}}} g_{k+1}(\beta_i^{(k)}).
\end{equation}

We can compute $\beta_0^{(k)},\beta_{k+1}^{(k)}$ by the method of Lagrange multipliers, which were first introduced in this general context by Pohst \cite{Pohst} (see Remark \ref{PohstLagrange}).  The values $a_{n-1},\dots,a_{n-k} \in \Z$ determine the power sums $s_i$ for $i=1,\dots,k$ by Newton's relations (\ref{Newton}).  Now the set of all $x=(x_i) \in \R^n$ such that $S_i(x)=s_i$ is closed and bounded, and therefore by symmetry the minimum (resp.\ maximum) value of the function $x_n$ on this set yields the bound $\beta_0^{(k)}$ (resp.\ $\beta_{k+1}^{(k)}$).  By the method of Lagrange multipliers, we find easily that if $x \in \R^n$ yields such an extremum, then there are at most $k-1$ distinct values among $x_1,\dots,x_{n-1}$, from which we obtain a finite set of possibilities for the extremum $x$.

For example, in the case $k=2$, the extrema are obtained from the equations
\[ (n-1)x_1+x_n=s_1=-a_{n-1} \quad\text{and}\quad (n-1)x_1^2+x_n^2=s_2=a_{n-1}^2-2a_{n-2} \]
which yields simply
\[ \beta_0^{(2)},\beta_3^{(2)}=\frac{1}{n}\left(-a_{n-1} \pm (n-1)\sqrt{a_{n-1}^2-2\left(1+\frac{1}{n-1}\right) a_{n-2}}\right). \]
(It is easy to show that this always improves upon the trivial bounds used by Takeuchi \cite{Takeuchi}.)  For $k=3$, for each partition of $n-1$ into $2$ parts, one obtains a system of equations which via elimination theory yield a (somewhat lengthy but explicitly given) degree $6$ equation for $x_n$.  For $k \geq 4$, we can continue in a similar way but we instead solve the system numerically, e.g., using the method of homotopy continuation as implemented by the package PHCpack developed by Verschelde \cite{Vers}; in practice, we do not significantly improve on these bounds whenever $k \geq 5$, and even for $k=5$, if $n$ is small then it often is more expensive to compute the improved bounds than to simply set $\beta_0^{(k)}=\beta_0^{(k-1)}$ and $\beta_{k+1}^{(k)}=\beta_k^{(k-1)}$.  

\begin{remark} \label{PohstLagrange}
Pohst's original use of Lagrange multipliers, which applies to number fields of arbitrary signature, instead sought the extrema of the power sum $S_{k+1}$ to bound the coefficient $a_{n-k-1}$.  The bounds given by Rolle's theorem for totally real fields are not only easier to compute (especially in higher degree) but in most cases turn out to be strictly stronger.  We similarly find that many other bounds typically employed in this situation (e.g., those arising from the positive definiteness of $T_2$ on $\Z[\alpha]$) are also always weaker.
\end{remark}

\subsection{Algorithmic details}

Our algorithm to solve Problem 1 then runs as follows.  We first apply the basic bounds from \S 2.2 to specify finitely many values of $a_{n-1},a_{n-2}$.  For each such pair, we use Rolle's theorem and the method of Lagrange multipliers to bound each of the coefficients inductively.  Note that if $k \geq 3$ is odd and $a_{n-1}=a_{n-3}=\dots=a_{n-(k-2)}=0$, then replacing $\alpha$ by $-\alpha$ we may assume that $a_{n-k} \geq 0$.

For each polynomial $f \in NS(n,B)$ that emerges from these bounds, we test it to see if it corresponds to a field $F \in NF(n,B)$.  We treat each of these latter two tasks in turn.

\subsubsection*{Calculation of real roots.}  

In the computation of the bounds (\ref{Rolle}), we use Newton's method to iteratively compute approximations to the roots $\beta_i^{(k)}$, using the fact that the roots of a polynomial are interlaced with those of its derivative, i.e. $\beta_{i-1}^{(k-1)} < \beta_i^{(k)}< \beta_i^{(k-1)}$ for $i=1,\dots,k$.  Note that by Rolle's theorem, we will either find a simple root in this open interval or we will converge to one of the endpoints, say $\beta_i^{(k)}=\beta_i^{(k-1)}$, and then necessarily $\beta_{i}^{(k)}=\beta_i^{(k-1)}$ as well, which implies that $f_k(x)$ is not squarefree and hence the entire coefficient range may be discarded immediately.  It is therefore possible to very quickly compute an approximate root which differs from the actual root $\beta_i^{(k+1)}$ by at most some fixed $\eps>0$.  We choose $\eps$ small enough to give a reasonable approximation but not so small as to waste time in Newton's method (say, $\eps=10^{-4}$).  We deal with the possibility of precision loss by bounding the value $g_{k+1}(\beta_i^{(k)})$ in (\ref{Rolle}) using elementary calculus; we leave the details to the reader.

\subsubsection*{Testing polynomials.}

For each $f \in NS(n,B)$, we test each of the following in turn.

\begin{enumerate}
\item We first employ an ``easy irreducibility test'': We rule out polynomials $f$ divisible by any of the factors: $x,x\pm 1,x\pm 2,x^2\pm x-1, x^2-2$.  In the latter three cases, we first evaluate the polynomial at an approximation to the values $(1 \pm \sqrt{5})/2,\sqrt{2}$, respectively, and then evaluate $f$ at these roots using exact arithmetic.  (Some benefit is gained by hard coding this latter evaluation.)  
\item We then compute the discriminant $d=\disc(f)$.  If $d \leq 0$, then $f$ is not a real separable polynomial, so we discard $f$.  
\item If $F=\Q[\alpha]=\Q[x]/(f(x)) \in NF(n,B)$, then for some $a \in \Z$ we have $B_O(n)^n < d_F=d/a^2 < B^n$ where $B_O$ is the Odlyzko bound (see \S 1).  Therefore using trial division we can quickly determine if there exists such an $a^2 \mid d$; if not, then we discard $f$.
\item Next, we check if $f$ is irreducible, and discard $f$ otherwise.
\item By the preceding two steps, an $a$-maximal order containing $\Z[\alpha]$ is in fact the maximal order $\Z_F$ of the field $F$.  If $\disc(\Z_F)=d_F > B$, we discard $f$.
\item Apply the \textsf{POLRED} algorithm of Cohen-Diaz y Diaz \cite{CDyDpol}: embed $\Z_F \subset \R^n$ by Minkowski (as in \S 1.1) and use LLL-reduction \cite{LLL} to compute a small element $\alpha_{\text{red}} \in \Z_F$ such that $\Q(\alpha)=\Q(\alpha_{\text{red}})=F$.  Add the minimal polynomial $f_{\text{red}}(x)$ of $\alpha_{\text{red}}$ to the list $NF(n,B)$ (along with the discriminant $d_F$), if it does not already appear.  
\end{enumerate}

We expect that almost all isomorphic fields will be identified in Step 6 by computing a reduced polynomial.  For reasons of efficiency, we wait until the space $NS(n,B)$ has been exhausted to do a final comparison with each pair of polynomials with the same discriminant to see if they are isomorphic.  Finally, we add the exceptional fields coming from Lemma \ref{Smyth}, if relevant.

\begin{remark}
Although Step 1 is seemingly trivial, it rules out a surprisingly significant number of polynomials $f$---indeed, nearly all reducible polynomials are discarded by this step in higher degrees.  Indeed, if $T_2(f)=\sum_i \alpha_i^2$ (where $\alpha_i$ are the roots of $f$) is small compared to $\deg(f)=n$, then $f$ is likely to be reducible and moreover divisible by a polynomial $g$ with $T_2(g)$ also small.  It would be interesting to give a precise statement which explains this phenomenon.
\end{remark}

\subsection{Implementation details}

For the implementation of our algorithm, we use the computer algebra system Sage \cite{SAGE}, which utilizes PARI \cite{Pari} for Steps 4--6 above.  Since speed was of the absolute essence, we found that the use of Cython (developed by Stein and Bradshaw) allowed us to develop a carefully optimized and low-level implementation of the bounds coming from Rolle's theorem and Lagrange multiplier method \ref{Rolle}.  We used the DSage package (due to Qiang) which allowed for the distribution of the compution to many machines; as a result, our computational time comes from a variety of processors (Opteron 1.8GHz, Athlon Dual Core 2.0GHz, and Celeron 2.53GHz), including a cluster of $30$ machines at the University of Vermont.

In low and intermediate degrees, where we expect comparatively many fields, we find that the running time is dominated by the computation of the maximal order (Step 5), followed by the check for irreducibility (Step 4); this explains the ordering of the steps as above.  By contrast, in higher degrees, where we expect few fields but must search in an exponentially large space, most of the time is spent in the calculation of real roots and in Step 1.  Further timing details can be found in Table 2 in \S 3.2.  

\section{Imprimitive fields}

In this section, we extend the ideas of the previous section to imprimitive fields $F$, i.e.\ those fields $F$ containing a nontrivial subfield.  Suppose that $F$ is an extension of $E$ with $[F:E]=m$ and $[E:\Q]=d$.  Since $\delta_F \geq \delta_E$, if $F \in NF(B)$ then $E \in NF(B)$ as well, and thus we proceed by induction on $E$.  For each such subfield $E$, we proceed in an analogous fashion.  We let 
\[ f(x)=x^m+a_{m-1}x^{m-1}+\dots+a_1 x + a_0 \] 
be the minimal polynomial of an element $\alpha \in \Z_F$ with $F=E(\alpha)$ and $a_i \in \Z_E$.

\subsection{Extension of bounds}

\subsubsection*{Basic bounds}

We begin with a relative version of Hunter's theorem.  We denote by $E_\infty$ the set of infinite places of $E$.  

\begin{lemma}[Martinet {\cite{MartinetRel}}] \label{relativeHunter}
There exists $\alpha \in \Z_F \setminus \Z_E$ such that
\begin{equation}
T_2(\alpha) \leq \frac{1}{m}\sum_{\sigma \in E_\infty} \left|\sigma\left(\Tr_{F/E} \alpha\right)\right|^2 + \gamma_{n-d} \left(\frac{|d_F|}{m^d |d_E|}\right)^{1/(n-d)}.
\label{relHunineq}
\end{equation}
\end{lemma}

The inequality of Lemma \ref{relativeHunter} remains true for any element of the set $\mu_E \alpha + \Z_E$, where $\mu_E$ denotes the roots of unity in $E$.  This allows us to choose $\Tr_{F/E} \alpha=-a_{m-1}$ among any choice of representatives from $\Z_E/m\Z_E$ (up to a root of unity); we choose the value of $a_{m-1}$ which minimizes 
\[ \sum_{\sigma \in E_\infty} \left|\sigma\left(\Tr_{F/E} \alpha\right)\right|^2=\sum_{\sigma \in E_\infty} \sigma(a_{m-1})^2, \] 
which is a positive definite quadratic form on $\Z_E$; such a value can be found easily using the LLL-algorithm.  

Now suppose that $F$ is totally real.  Then $\sum_{\sigma \in E_\infty} \left|\sigma(a_{m-1})\right|^2=\Tr_{E/\Q} a_{m-1}^2$, and we have $2^d$ or $\lceil m^d/2 \rceil$ possibilities for $a_{m-1}$, according as $m=2$ or otherwise.  For each value of $a_{m-1}$, we have $T_2(\alpha) \in \Z$ bounded from above by Lemma \ref{relativeHunter} and from below by Lemma \ref{Smyth} since $\Tr_{E/\Q}(\alpha^2)=T_2(\alpha) > 1.7719n$.  If we denote $\Tr_{F/E} \alpha^2=t_2$, then by Newton's relations, we have $t_2=a_{m-1}^2-2a_{m-2}$, and hence $\Tr_{E/\Q} t_2=T_2(\alpha)$ and $t_2 \equiv a_{m-1}^2 \pmod{2}$.  In particular, $t_2 \in \Z_E$ is totally positive and has bounded trace, leaving only finitely many possibilities: indeed, if we embed $\Z_E \hookrightarrow \R^d$ by Minkowski, these inequalities define a parallelopiped in the positive orthant.  

\subsubsection*{Lattice points in boxes.}

One option to enumerate the possible values of $t_2$ is to enumerate all lattice points in a sphere of radius given by (\ref{relHunineq}) using the Fincke-Pohst algorithm \cite{FinckePohst}.  However, one ends up enumerating far more than what one needs in this fashion, and so we look to do better.  The problem we need to solve is the following.

\begin{problem} \label{lattic}
Given a lattice $L \subset \R^d$ of rank $d$ and a convex polytope $P$ of finite volume, enumerate the set $P \cap L$.  
\end{problem}

Here we must allow the lattice $L$ to be represented numerically; to avoid issues of precision loss, one supposes without loss of generality that $\partial P \cap L = \emptyset$.

There exists a vast literature on the classical problem of the enumeration of integer lattice points in rational convex polytopes (see e.g., De Loera \cite{DeLoera}), as well as several implementations \cite{LattE,PALP}.  (In many cases, these authors are concerned primarily with simply counting the number of lattice points, but their methods equally allow their enumeration.)

In order to take advantage of these methods to solve Problem \ref{lattic}, we compute an LLL-reduced basis $\gamma=\gamma_1,\dots,\gamma_d$ of $L$, and we perform the change of variables $\phi: \R^d \to \R^d$ which maps $\gamma_i \mapsto e_i$ where $e_i$ is the $i$th coordinate vector.  The image $\phi(P)$ is again a convex polytope.  We then compute a rational polytope $Q$ (i.e. a polytope with integer vertices) containing $\phi(P)$ by rounding the vertices to the nearest integer point as follows.  For each pair of vertices $v,w \in P$ such that the line $\ell(v,w)$ containing $v$ and $w$ is not contained in a proper face of $P$, we round the $i$th coordinates $\phi(v)_i$ down and $\phi(w)_i$ up if $\phi(v)_i \leq \phi(w)_i$, and otherwise round in the opposite directions.  The convex hull $Q$ of these rounded vertices clearly contains $\phi(P \cap L)$, and is therefore amenable to enumeration using the methods above.

We note that in the case where $P$ is a parallelopiped, for each vertex $v$ there is a unique opposite vertex $w$ such that the line $\ell(v,w)$ is not contained in a proper face, so the convex hull $Q$ will also form a parallelopiped.  

\subsubsection*{Coefficient bounds and testing polynomials.}

The bounds in \S 2 apply \emph{mutatis mutandis} to the relative situation.  For example, given $a_{m-1},\dots,a_{m-k}$ for $k \geq 2$, for each $v \in E_\infty$, if we let $v(g)$ denote the polynomial $\sum_i v(b_i) x^i$ for $g(x)=\sum_i b_i x^i \in E[x]$, we obtain the inequality
\[ -\min_{\substack{0 \leq i \leq k+1 \\ i \not\equiv k \psod{2}}} v(g_{k+1})(\beta_{i,v}^{(k)}) < v(a_{m-k-1}) < - \max_{\substack{0 \leq i \leq k+1 \\ i \equiv k \psod{2}}} v(g_{k+1})(\beta_{i,v}^{(k)}); \]
here, $\beta_{1,v}^{(k)},\dots,\beta_{k,v}^{(k)}$ denote the roots of $v(f_k(x))$, and $\beta_{0,v}^{(k)},\beta_{k+1,v}^{(k)}$ are computed in an analogous way using Lagrange multipliers.  In this situation, we have $a_{m-k-1}$ contained in an honest rectangular box, and the results of the previous subsection apply directly.

For each polynomial which satisfies these bounds, we perform similar tests to discard polynomials as in \S 2.3.  One has the option of working always relative to the ground field or immediately computing the corresponding absolute field; in practice, for the small base fields under consideration, these approaches seem to be comparable, with a slight advantage to working with the absolute field.

\subsection{Conclusion and timing}

Putting together the primitive and imprimitive fields computed in \S\S 2--3, we have proven Theorem 2.  In Table 2, we list some timing details arising from the computation.  Note that in high degrees (presumably because we enumerate an exponentially large space) we recover all imprimitive fields already during the search for primitive fields.

\begin{center} \label{table:timing}
\begin{table}
\caption{Timing data}
\begin{tabular}{c|ccccccccc}
$n$                          & 2    & 3    & 4        & 5        & 6         & 7        & 8         & 9         & 10   \\ 
\hline
$\Delta(n)$                  & 30   & 25   & 20       & 17       & 16        & 15.5     & 15        & 14.5      & 14   \\
$f$                          & 443  & 4922 & 57721    & 244600   & 3242209   & $1.7 \times 10^7$ & $1.2\times 10^8$ & $9.5\times 10^8$ & $2.2\times 10^9$ \\
Irred $f$                    & 418  & 2523 & 27234    & 157613   & 2710965   & $1.6 \times 10^7$ & $1.1\times 10^8$ & $9.0\times 10^8$ & $2.2\times 10^9$ \\
$f$, $d_F \leq B$            & 418  & 1573 & 5665     & 4497     & 1288      & 4839     & 3016      & 506       & 0 \\
$F$                          & 273  & 630  & 1273     & 674      & 802       & 301      & 164       & 15        & 0 \\
Total time                   & 0.2s & 2.2s & 26.8s    & 1m25s    & 17m3s     & 2h59m    & 1d4.5h    & 17d21h    & 173d \\
\hline
Imprim $f$                   & 0    & 0    & 7059     & 0        & 62532     & 0        & 239404    & 15658     & 945866 \\
Imprim $F$                   & 0    & 0    & 702      & 0        & 420       & 0        & 100       & 6         & 0 \\
Time                         & -    & -    & 4m22s    & -        & 8m38s     & -        & 1h56m     & 16m53s    & 11h27m \\
\hline
Total fields                 & 273  & 630  & 1578     & 674      & 827       & 301      & 164       & 15        & 0
\end{tabular}
\end{table}
\end{center}

\section{Tables of totally real fields}

In Table 3, we count the number of totally real fields $F$ with root discriminant $\delta_F \leq 14$ by degree, and separate out the primitive and imprimitive fields.  We also list the minimal discriminant and root discriminant for $n \leq 9$.  The polynomial
\[ x^{10} - 11x^8-3x^7+37x^6+14x^5-48x^4-22x^3+20x^2+12x+1 \]
with $d_F=443952558373=61^2 397^2 757$ and $\delta_F \approx 14.613$ is the dectic totally real field with smallest discriminant that we found---the corresponding number field (though not this polynomial) already appears in the tables of Kl\"uners-Malle \cite{KM} and is a quadratic extension of the second smallest real quintic field, of discriminant $24217$.  It is reasonable to conjecture that this is indeed the smallest such field.

\begin{table}[h] \label{table:data}
\caption{Totally real fields $F$ with $\delta_F \leq 14$}
\begin{center}
\begin{tabular}{c|c|cc|cc}
$n=[F:\Q]$ &\ $\#NF(n,14)$\ \ &\ Primitive $F$\ &\ Imprimitive $F$\ \ & Minimal $d_F$ & Minimal $\delta_F$ \\
\hline
\rule{0pt}{2.5ex} 2 \rule{0pt}{2.5ex}  & 59 & 59 & 0 & 5 & 2.236 \\
3 & 86 & 86 & 0 & 49 & 3.659 \\
4 & 277 & 117 & 160 & 725 & 5.189 \\
5 & 170 & 170 & 0 & 14641 & 6.809 \\
6 & 263 & 104 & 159 & 300125 & 8.182 \\
7 & 301 & 301 & 0 & 20134393 & 11.051 \\
8 & 62 & 19 & 43 & 282300416 & 11.385 \\
9 & 11 & 6 & 5 & 9685993193 & 12.869 \\
10 & 0 & 0 & 0 & \ 443952558373? & 14.613? \\
\hline
\rule{0pt}{2.5ex} 
Total \rule{0pt}{2.5ex}  & 1229 & 862 & 367 & - & -
\end{tabular}
\end{center}
\end{table}

In Tables 4--5, we list the octic and nonic fields $F$ with $\delta_F \leq 14$.  For each field, we specify a maximal subfield $E$ by its discriminant and degree---when more than one such subfield exists, we choose the one with smallest discriminant.  

\begin{table} \label{table:octic1}
\caption{Octic totally real fields $F$ with $\delta_F \leq 14$}
\begin{scriptsize}
\[
\begin{array}{cccc}
d_F & f & [E:\Q] & d_E \\ \hline
282300416 & x^8 - 4x^7 + 14x^5 - 8x^4 - 12x^3 + 7x^2 + 2x - 1 & 4 & 2624 \\
309593125 & x^8 - 4x^7 - x^6 + 17x^5 - 5x^4 - 23x^3 + 6x^2 + 9x - 1 & 4 & 725 \\
324000000 & x^8 - 7x^6 + 14x^4 - 8x^2 + 1 & 4 & 1125 \\
410338673 & x^8 - x^7 - 7x^6 + 6x^5 + 15x^4 - 10x^3 - 10x^2 + 4x + 1 & 4 & 4913 \\
432640000 & x^8 - 2x^7 - 7x^6 + 16x^5 + 4x^4 - 18x^3 + 2x^2 + 4x - 1 & 4 & 1600 \\
442050625 & x^8 - 2x^7 - 12x^6 + 26x^5 + 17x^4 - 36x^3 - 5x^2 + 11x - 1 & 4 & 725 \\
456768125 & x^8 - 2x^7 - 7x^6 + 11x^5 + 14x^4 - 18x^3 - 8x^2 + 9x - 1 & 4 & 725 \\
483345053 & x^8 - x^7 - 7x^6 + 4x^5 + 15x^4 - 3x^3 - 9x^2 + 1 & 1 & 1 \\
494613125 & x^8 - x^7 - 7x^6 + 4x^5 + 13x^4 - 4x^3 - 7x^2 + x + 1 & 4 & 725 \\
582918125 & x^8 - 2x^7 - 6x^6 + 9x^5 + 11x^4 - 9x^3 - 6x^2 + 2x + 1 & 4 & 725 \\
656505625 & x^8 - 3x^7 - 4x^6 + 13x^5 + 5x^4 - 13x^3 - 4x^2 + 3x + 1 & 4 & 725 \\
661518125 & x^8 - x^7 - 7x^6 + 5x^5 + 15x^4 - 7x^3 - 10x^2 + 2x + 1 & 2 & 5 \\
707295133 & x^8 - 8x^6 - 2x^5 + 19x^4 + 7x^3 - 13x^2 - 4x + 1 & 1 & 1 \\
733968125 & x^8 - 2x^7 - 6x^6 + 10x^5 + 11x^4 - 11x^3 - 7x^2 + 2x + 1 & 2 & 5 \\
740605625 & x^8 - x^7 - 9x^6 + 8x^5 + 21x^4 - 12x^3 - 14x^2 + 4x + 1 & 4 & 725 \\
803680625 & x^8 - 2x^7 - 9x^6 + 12x^5 + 22x^4 - 24x^3 - 14x^2 + 14x - 1 & 4 & 725 \\
852038125 & x^8 - 10x^6 - 5x^5 + 17x^4 + 5x^3 - 10x^2 + 1 & 4 & 725 \\
877268125 & x^8 - 3x^7 - 6x^6 + 20x^5 + 5x^4 - 25x^3 - x^2 + 7x + 1 & 4 & 725 \\
898293125 & x^8 - x^7 - 9x^6 + 10x^5 + 15x^4 - 10x^3 - 9x^2 + x + 1 & 4 & 725 \\
1000118125 & x^8 - 3x^7 - 4x^6 + 14x^5 + 5x^4 - 19x^3 - x^2 + 7x - 1 & 2 & 5 \\
1024000000 & x^8 - 8x^6 + 19x^4 - 12x^2 + 1 & 4 & 1600 \\
1032588125 & x^8 - 9x^6 - 2x^5 + 23x^4 + 9x^3 - 17x^2 - 9x - 1 & 2 & 5 \\
1064390625 & x^8 - 13x^6 + 44x^4 - 17x^2 + 1 & 4 & 725 \\
1077044573 & x^8 - x^7 - 8x^6 + 8x^5 + 16x^4 - 17x^3 - 2x^2 + 5x - 1 & 1 & 1 \\
1095205625 & x^8 - 3x^7 - 5x^6 + 18x^5 + 2x^4 - 23x^3 + 2x^2 + 8x - 1 & 2 & 5 \\
1098290293 & x^8 - 3x^7 - 4x^6 + 16x^5 + x^4 - 23x^3 + 7x^2 + 5x - 1 & 1 & 1 \\
1104338125 & x^8 - 2x^7 - 8x^6 + 15x^5 + 17x^4 - 31x^3 - 9x^2 + 17x - 1 & 4 & 725 \\
1114390153 & x^8 - 8x^6 - 2x^5 + 16x^4 + 3x^3 - 10x^2 + 1 & 1 & 1 \\
1121463125 & x^8 - 3x^7 - 4x^6 + 15x^5 + 2x^4 - 18x^3 + 5x + 1 & 2 & 5 \\
1136700613 & x^8 - x^7 - 7x^6 + 4x^5 + 14x^4 - 4x^3 - 8x^2 + x + 1 & 1 & 1 \\
1142440000 & x^8 - 3x^7 - 5x^6 + 15x^5 + 8x^4 - 15x^3 - 5x^2 + 4x + 1 & 4 & 4225 \\
1152784549 & x^8 - 4x^7 - x^6 + 15x^5 - 3x^4 - 16x^3 + 4x^2 + 4x - 1 & 4 & 1957 \\
1153988125 & x^8 - 2x^7 - 7x^6 + 11x^5 + 12x^4 - 16x^3 - 5x^2 + 6x - 1 & 4 & 2525 \\
1166547493 & x^8 - x^7 - 7x^6 + 6x^5 + 14x^4 - 9x^3 - 9x^2 + 3x + 1 & 1 & 1 \\
1183423341 & x^8 - x^7 - 8x^6 + 9x^5 + 17x^4 - 20x^3 - 8x^2 + 10x - 1 & 4 & 1957 \\
1202043125 & x^8 - 3x^7 - 4x^6 + 16x^5 - 21x^3 + 9x^2 + 2x - 1 & 2 & 5 \\
1225026133 & x^8 - 3x^7 - 4x^6 + 18x^5 - 6x^4 - 17x^3 + 9x^2 + 2x - 1 & 1 & 1 \\
1243893125 & x^8 - x^7 - 8x^6 + 3x^5 + 18x^4 - x^3 - 12x^2 - 2x + 1 & 2 & 5 \\
1255718125 & x^8 - 2x^7 - 8x^6 + 19x^5 + 10x^4 - 41x^3 + 13x^2 + 10x - 1 & 4 & 725 \\
1261609229 & x^8 - 2x^7 - 6x^6 + 12x^5 + 9x^4 - 19x^3 - x^2 + 6x - 1 & 1 & 1 \\
1292203125 & x^8 - 4x^7 - x^6 + 17x^5 - 6x^4 - 21x^3 + 6x^2 + 8x + 1 & 4 & 1125 \\
1299600812 & x^8 - 2x^7 - 6x^6 + 10x^5 + 12x^4 - 13x^3 - 8x^2 + 3x + 1 & 1 & 1 \\
1317743125 & x^8 - x^7 - 8x^6 + 7x^5 + 19x^4 - 14x^3 - 12x^2 + 8x - 1 & 2 & 5 \\
1318279381 & x^8 - x^7 - 7x^6 + 5x^5 + 14x^4 - 6x^3 - 9x^2 + x + 1 & 1 & 1 \\
1326417388 & x^8 - 2x^7 - 6x^6 + 10x^5 + 12x^4 - 13x^3 - 9x^2 + 4x + 2 & 4 & 2777 \\
1348097653 & x^8 - 2x^7 - 6x^6 + 11x^5 + 11x^4 - 17x^3 - 6x^2 + 6x + 1 & 1 & 1 \\
1358954496 & x^8 - 8x^6 + 20x^4 - 16x^2 + 1 & 4 & 2048 \\
1359341129 & x^8 - 8x^6 - x^5 + 18x^4 + 2x^3 - 12x^2 - x + 2 & 1 & 1 \\
1377663125 & x^8 - 12x^6 + 33x^4 - 5x^3 - 22x^2 + 5x + 1 & 4 & 725 \\
1381875749 & x^8 - 3x^7 - 4x^6 + 14x^5 + 4x^4 - 18x^3 + x^2 + 5x - 1 & 1 & 1 \\
1391339501 & x^8 - 3x^7 - 4x^6 + 15x^5 + 4x^4 - 22x^3 + 9x - 1 & 1 & 1 \\
1405817381 & x^8 - 9x^6 - x^5 + 20x^4 + 6x^3 - 12x^2 - 7x - 1 & 1 & 1 \\
1410504129 & x^8 - 9x^6 - x^5 + 22x^4 + x^3 - 15x^2 - x + 1 & 4 & 3981 \\
1410894053 & x^8 - 2x^7 - 6x^6 + 9x^5 + 12x^4 - 11x^3 - 8x^2 + 3x + 1 & 1 & 1 \\
1413480448 & x^8 - 4x^7 - 2x^6 + 16x^5 - x^4 - 16x^3 + 2x^2 + 4x - 1 & 4 & 2048 \\
1424875717 & x^8 - x^7 - 7x^6 + 5x^5 + 15x^4 - 6x^3 - 10x^2 + x + 1 & 1 & 1 \\
1442599461 & x^8 - 3x^7 - 4x^6 + 15x^5 + 4x^4 - 21x^3 - 2x^2 + 8x + 1 & 4 & 7053 \\
1449693125 & x^8 - x^7 - 9x^6 + 10x^5 + 20x^4 - 20x^3 - 14x^2 + 11x + 1 & 2 & 5 \\
1459172469 & x^8 - 4x^7 - x^6 + 17x^5 - 6x^4 - 21x^3 + 8x^2 + 6x - 1 & 4 & 1957 \\
1460018125 & x^8 - 3x^7 - 5x^6 + 13x^5 + 11x^4 - 14x^3 - 10x^2 + x + 1 & 4 & 2525 \\
1462785589 & x^8 - 2x^7 - 6x^6 + 11x^5 + 10x^4 - 17x^3 - 3x^2 + 6x - 1 & 1 & 1 \\
1472275625 & x^8 - 3x^7 - 6x^6 + 19x^5 + 13x^4 - 35x^3 - 12x^2 + 13x - 1 & 4 & 725 \\
\end{array}
\]
\end{scriptsize}
\end{table}

\begin{scriptsize}
\begin{table}[ht] \label{table:nonic}
\caption{Nonic totally real fields $F$ with $\delta_F \leq 14$}
\[ 
\begin{array}{cccc}
d_F & f & [E:\Q] & d_E \\ \hline
\noalign{\smallskip}
9685993193 & x^9 - 9x^7 + 24x^5 - 2x^4 - 20x^3 + 3x^2 + 5x - 1 & 1 & 1 \\
11779563529 & x^9 - 9x^7 - 2x^6 + 22x^5 + 5x^4 - 17x^3 - 4x^2 + 4x + 1 & 1 & 1 \\
16240385609 & x^9 - x^8 - 9x^7 + 4x^6 + 26x^5 - 2x^4 - 25x^3 - x^2 + 7x + 1 & 3 & 49 \\
16440305941 & x^9 - 2x^8 - 9x^7 + 11x^6 + 28x^5 - 18x^4 - 34x^3 + 8x^2 + 13x + 1 & 3 & 229 \\
16898785417 & x^9 - 2x^8 - 7x^7 + 11x^6 + 18x^5 - 17x^4 - 19x^3 + 6x^2 + 7x + 1 & 1 & 1 \\
16983563041 & x^9 - x^8 - 8x^7 + 7x^6 + 21x^5 - 15x^4 - 20x^3 + 10x^2 + 5x - 1 & 3 & 361 \\
17515230173 & x^9 - 4x^8 - 3x^7 + 29x^6 - 26x^5 - 24x^4 + 34x^3 - 2x^2 - 5x + 1 & 3 & 49 \\
18625670317 & x^9 - 9x^7 - x^6 + 23x^5 + 4x^4 - 19x^3 - 3x^2 + 4x + 1 & 1 & 1 \\
18756753353 & x^9 - 3x^8 - 4x^7 + 15x^6 + 4x^5 - 22x^4 - x^3 + 10x^2 - 1 & 1 & 1\\
19936446593 & x^9 - 3x^8 - 5x^7 + 17x^6 + 7x^5 - 30x^4 - x^3 + 16x^2 - 2x - 1 & 3 & 49 \\
20370652633 & x^9 - 2x^8 - 8x^7 + 12x^6 + 15x^5 - 17x^4 - 8x^3 + 8x^2 + x - 1 & 1 & 1\\
\end{array}
\]
\end{table}
\end{scriptsize}

\end{document}